\documentclass[12pt, epsf]{article}
\usepackage{latexsym}
\usepackage{amsmath}
\usepackage{amssymb}
\usepackage{amsfonts}
\usepackage{times}
\usepackage{graphicx}
\usepackage{epsfig}
\usepackage{array}
\usepackage{flafter}
\usepackage{nicefrac}
\newtheorem{definition}{Definition}
\newtheorem{theorem}{Theorem}
\newtheorem{ExampleDef}{Example}[section]

\newcommand{\Example}[3]{
  \begin{list}{}{
      \setlength{\leftmargin}{1em}} % Indent everything by this amount
    \item                           % Group everything in one item
    \small                          % Use a smaller font size
    \begin{ExampleDef} \rm          % select roman
      {\bf \hspace{-1ex}: #1}       % Use[1ex] to break line
      #2                                % The actual stuff
      \hfill {\large \boldmath $\Box$}  % The box
      \label{ex#3}                      % Label the example
    \end{ExampleDef}
  \end{list}}

\begin{document}

\begin{center}
{\Large \bf Notes on stable learning\\ with piecewise-linear basis functions \par} \vspace{1.5em}
{\large Winfried Lohmiller, Philipp Gassert, and Jean-Jacques Slotine \par}
{Nonlinear Systems Laboratory \\
Massachusetts Institute of Technology \\
Cambridge, Massachusetts, 02139, USA\\
{\sl wslohmil@mit.edu} \par}
\vspace{3em}
\end{center}

\section{Contraction theory} \label{contraction}

Our main stability analysis tool is nonlinear contraction theory~\cite{Lohm1}.
Specifically, we will use this specialized result.

\begin{theorem}
\itshape Consider the system $\dot{\bf x} = {\bf f(x},t)$, where ${\bf
  f}$ is a differentiable nonlinear function of ${\bf x} \in R^n$ and time $t$,
and a constant invariant set for this dynamics.  Assume that, in this
set, the system is contracting in an identity metric, i.e. that its
Jacobian is uniformly negative definite. Then all trajectories starting
in the invariant set converge exponentially to each other, with a rate
equal to the positive constant $\lambda$, where $\ - \lambda {\bf I}
\ $ is an upper bound on the Jacobian.
 \label{th:theoremF}
\end{theorem}

When the Jacobian is only negative semi-definite, the system is said to be
{\it semi-contracting}.  Then the length of any path between two trajectories within the invariant set can only decrease.

\section{Learning}  

Over the last decades much progress has been achieved in the field of
neuro-inspired learning, and in particular deep
learning~\cite{Hinton}. In these algorithms, the measurement equation
is typically based on a piece-wise continuous function
\begin{equation}
y_k = \sum_{i=1}^I f_i \left( {\bf w}_i^T {\bf x}_k  \right) \nonumber
\end{equation}
with the constant unknown parameters ${\bf w}_i = (w_{i1}, ..., w_{iN+1})^T$ of the $i = 1..I$ neurons, the $k=1..K$ inputs ${\bf x}_k(t) = (1, x_{1k}^T, .., x_{Nk})^T$ and outputs $y_k(t)$, and the $n = 1..N$ input elements $x_{nk}$. Note that the first element of ${\bf x}_k(t)$, which is always $1$, corresponds to a constant offset. We assume now and in the following that $f_i$ is piece-wise linear.

We want to estimate this function with
\begin{equation}
\hat{y}_k = \sum_{i=1}^I f_i \left( \hat{z}_{ik}  \right) \nonumber 
\end{equation}
and $\hat{z}_{ik} =  \hat{\bf w}_i^T {\bf x}_k$. The standard way to estimate $\hat{\bf w}_i$ is to define and minimize the cost function
\begin{equation}
V =  \frac{1}{2} \sum_ {k=1}^K \tilde{y}_k^2 \label{eq:V}
\end{equation}
with $\tilde{y}_k  = \hat{y}_k - y_k$.   The parameter learning law is in the direction of the gradient of (\ref{eq:V})
\begin{equation}
{\dot{\hat{\bf w}}}_i = - \frac{1}{T}  \frac{\partial V}{\partial \hat{\bf w}_i} = - \frac{1}{T} \sum_ {k=1}^K \tilde{y}_k \frac{\partial \hat{y}_k}{\partial \hat{z}_{ik}} {\bf x}_k \label{eq:dotb}
\end{equation}
with $T(t) > 0$
This learning law is a set of switched linear functions for piece-wise linear $f_i$. Note that for $K=1$, constant ${\bf x}$ and $T = {\bf x}^T {\bf x}$ the above can be simplified to
\begin{equation}
{\dot{\hat{z}}}_{i} = -  \tilde{y} \frac{\partial \hat{y}}{\partial \hat{z}_{i}}  \label{eq:dotz}
\end{equation}
where we have omited the index $k$. The Jacobian of the above is the block diagonal matrix
\begin{equation}
\frac{\partial \dot{\hat{\bf w}}_i}{\partial \hat{\bf w}_j} = - \frac{1}{T} \  \sum_ {k=1}^K {\bf x}_k {\bf x}_k^T \frac{\partial^2 V_k}{\partial \hat{\bf z}_{ik} \partial \hat{\bf z}_{jk}} \nonumber 
\end{equation}
with the Hessian 
\begin{equation}
\frac{\partial^2 V_k}{\partial \hat{\bf z}_{ik} \partial \hat{\bf z}_{jk}}  = \frac{\partial \hat{y}_{k}}{\partial \hat{\bf z}_{ik}} \frac{\partial \hat{y}_k}{\partial \hat{\bf z}_{ik}}^T + \frac{\partial^2 \hat{y}_k }{\partial \hat{\bf z}_{ik} \partial \hat{\bf z}_{jk}} \tilde{y}_k \nonumber
\end{equation}
The first dyadic term is n.s.d. which implies convexity of the cost function in its local linear region. Each linear region is semi-contracting and contracting in the observability direction $\frac{\partial \hat{y}_k}{\partial \hat{\bf z}_{k}}  {\bf x}_k$ 

The second term vanishes everywhere except at the transition of one linear region to its neighboring linear region. At this transition the definiteness of this term changes its sign with the sign of $\tilde{y}$ even if  the diagonal Hessian $\frac{\partial^2 \hat{y}_k }{\partial \hat{\bf z}_{ik} \partial \hat{\bf z}_{jk}}$ is semi-definite. Note that according to Theorem \ref{th:theoremF} local diverging areas are unproblematic, as long as they do not separate a finite region from the contracting space.

Section \ref{standard} explains that for standard saturation functions $f_i$ finite regions of the state space are separated by diverging transitions to the main semi-contraction regions. Hence finite regions of the state space will not converge to the true solution. This problem is equivalent to local minima and shallow areas of the cost function and can lead to failures of learning, see e.g. \cite{Shai}. Both local minima and shallow areas are caused by the non-convexity of the cost function.

\begin{figure}
\begin{center}
\includegraphics[scale=0.25]{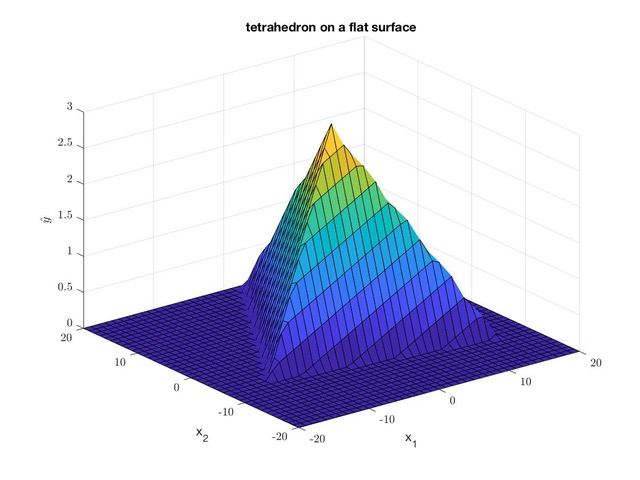}
\end{center}
\caption{Tetrahedron on flat plane}
\label{fig:pyr}
\end{figure}

In addition are neuronal networks in general strongly under-determined, i.e. there are much more neurons than degrees of freedom of the observed plant. This increases computation time and implies the problem of learning non-real effects due to noise. In Kalman theory such under-determined observers are not recommended due to that reason. Also, standard neurons like ReLUs (Rectified Linear Unit), as $\max(0, \hat{z})$, always generate an infinitesimal edge at $ \hat{z} = 0$. Hence they cannot approximate general piece-wise linear surfaces with finite edges as e.g. the pyramid in figure \ref{fig:pyr}.

To resolve this problem section \ref{nn} introduces a function approximation, which does not separate finite regions of the state space to the semi-contraction region. In the proposed new approach we try to estimate a non-linear function with the minimal set of neurons to reduce this complexity. This however requires that every neuron converges to its right weights, even if only few neurons exist. 

Also, we have to select piece-wise linear neurons with finite edges to be able to approximate any piece-wise linear function.  Section \ref{jn} will do so for $N$-dimensional surfaces where the line segments of $1$-dimensional functions turn to $N$-dimensional hyperplane segements.

\subsection{Standard learning of $1$-dimensional functions} \label{standard}

We illustrate the problem of learning with saturation functions (ReLUs) for the special case of 2 positive neurons $f_1(\hat{z}_1) = \max(0, \hat{z}_1)$,  $f_2(\hat{z}_2) = \max(0, \hat{z}_2)$ with $\hat{z}_1 = \hat{w}_{11} + x \hat{w}_{12}$, $\hat{z}_2 = \hat{w}_{21} + x \hat{w}_{22}$  
which corresponds to the estimated function
\begin{equation}
\hat{y} = 
\left\{
\begin{array}{l}
\hat{z}_1 + \hat{z}_2    \mathrm{:for } \ \hat{z}_1 \ge 0, \hat{z}_2  \ge 0  \\
\hat{z}_2     \mathrm{:for } \ \hat{z}_1 <0, \hat{z}_2 > 0 \\
\hat{z}_1    \mathrm{:for } \ \hat{z}_1 >0, \hat{z}_2 < 0  \\
0  \mathrm{:for } \ \hat{z}_1 \le 0, \hat{z}_2  \le 0
\end{array} \right. \nonumber
\end{equation}
with 
\begin{equation}
\frac{\partial \hat{y}}{\partial \hat{\bf z} } = 
\left\{
\begin{array}{l}
(1, 1)^T    \mathrm{:for  } \ \hat{z}_1 \ge 0, \hat{z}_2  \ge 0  \\
(0, 1)^T     \mathrm{:for } \ \hat{z}_1 <0, \hat{z}_2 > 0 \\
(1, 0)^T    \mathrm{:for  } \ \hat{z}_1 >0, \hat{z}_2 < 0  \\
(0, 0)^T  \mathrm{:for } \ \hat{z}_1 \le 0, \hat{z}_2  \le 0
\end{array} \right. \nonumber
\end{equation}
At the transitions between the linear regions we get the semi-definite Jacobian
\begin{equation}
\frac{\partial^2 \hat{y}}{\partial \hat{\bf z} \partial \hat{\bf z}}  = 
\left\{
\begin{array}{l}
\left(
\begin{array}{ccc}
\delta  &  0 \\
0  & 0
\end{array}
\right) \mathrm{:for } \ \hat{z}_1 = 0  \\
\left(
\begin{array}{ccc}
 0 &  0 \\
0  & \delta
\end{array}
\right)    \mathrm{:for } \ \hat{z}_2 = 0 \\
0  \mathrm{:otherwise} 
\end{array} \right. \label{eq:transitionsat}
\end{equation}
with the Dirac impulse $\delta > 0$. Note these definite matrices are then weighted with $\tilde{y}$, i.e. they are semi-contracting (semi-diverging) for $\tilde{y}>0$ ($\tilde{y}<0$).

%\begin{figure}
%\begin{center}
%\includegraphics[scale=0.5]{cc2p0.jpg}
%\end{center}
%\caption{Flow field of 2 saturation neurons and $y=0$}
%\label{fig:t1}
%\end{figure}

\begin{figure}
\begin{center}
\includegraphics[scale=0.3]{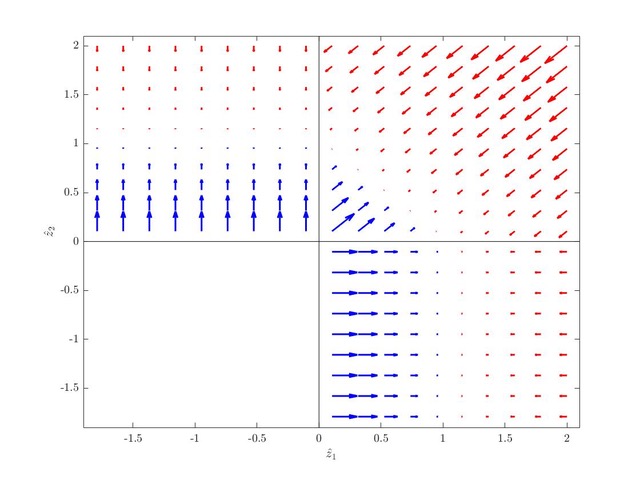}
\end{center}
\caption{Flow field of 2 saturation neurons and $y=1$}
\label{fig:t2}
\end{figure}

\begin{figure}
\begin{center}
\includegraphics[scale=0.3]{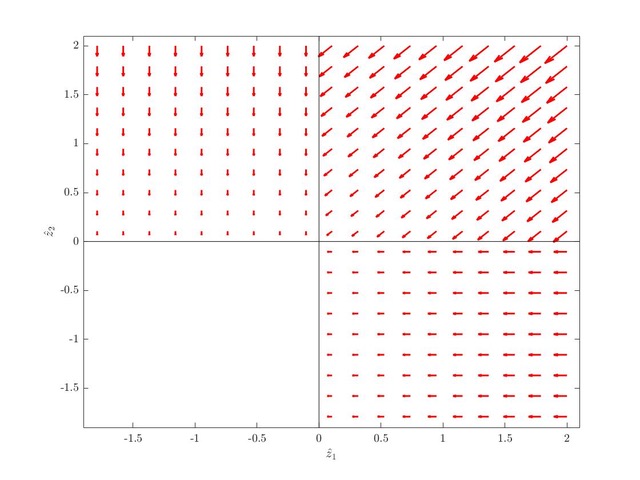}
\end{center}
\caption{Flow field of 2 saturation neurons  and $y=-1$}
\label{fig:t3}
\end{figure}

We illustrate this now for the simplified learning law (\ref{eq:dotz}) with 2 saturation neurons in figure \ref{fig:t2} and \ref{fig:t3} respectively. Red coloring means that $\tilde{y}>0$ and blue coloring means that $\tilde{y}<0$.

Figure \ref{fig:t2} shows that the lower left corner is separated by diverging transitions  (\ref{eq:transitionsat}) for blue $\tilde{y}<0$. Figure \ref{fig:t3} has only semi-contracting transitions (\ref{eq:transitionsat}) for red $\tilde{y}>0$. However all trajectories are moved to the lower left corner, where they freeze then indifferently.  In both cases the major problem is that there are separated finite regions to the semi-contraction region of red $\tilde{y}>0$. Hence the learning can be stuck at an incorrect weight. 

At this example the separated regions can only be avoided if at most one transition between linear regions exists. Hence the following sections proposes a learning law with at most one transition per neuron.

\subsection{Contracting learning of 1-dimensional functions} \label{nn}

In the following we propose a 1-dimensional approximation function which is globally (semi)-contracting except for regions of zero measure. Indeed we locally estimate the piece-wise linear function 
\begin{equation}
\hat{y}_k =  \hat{z}_{ik}  \  \mathrm{: for \  } x_{i-1} \le x \le  x_{i}    \label{eq:lin} \\ 
\end{equation}
with the local linear neuron $\hat{z}_{ik} =  \hat{\bf w}_i^T {\bf x}_k$ with $\hat{\bf w}_i = (\hat{w}_{i1}, \hat{w}_{i2})^T$ and ${\bf x}_k = (1, x_k)^T$ between the  corners $x_{i-1}, x_{i}$ to the neighboring linear neurons $i-1$ and $i+1$, that solve $\hat{z}_{i} = \hat{z}_{i-1}, \hat{z}_{i} = \hat{z}_{i+1}$ respectively. 

We can now compute
\begin{equation}
\frac{\partial \hat{y}_k}{\partial \hat{z}_{ik}} = 
\left\{
\begin{array}{l}
1  \mathrm{:for \ active \ neurons \ }  \\
0  \mathrm{:for \ non \ active \ neurons \ } 
\end{array} \right. \nonumber
\end{equation}
At the transitions between the linear regions we get the semi-definite Jacobian
\begin{equation}
\frac{\partial^2 \hat{y}_k}{\partial \hat{\bf z}_{ik} \partial \hat{\bf z}_{jk}}  = 
\left\{
\begin{array}{l}
\left(
\begin{array}{ccc}
\delta  &  - \delta  \\
- \delta   &  \delta 
\end{array}
\right) \mathrm{: for \ the  \ corner \ between \ neuron \ } i \ \mathrm{to \ its  \ neighbour } \\
0   \mathrm{: otherwise}  \\
\end{array} \right. \nonumber
\end{equation}
with the Dirac impulse $\delta > 0$. Note these definite matrices are then weighted with $\tilde{y}$, i.e. they are semi-contracting (semi-diverging) for $\tilde{y}>0$ ($\tilde{y}<0$).

%\begin{figure}
%\begin{center}
%\includegraphics[scale=0.5]{ccadd0.jpg}
%\end{center}
%\caption{Flow field of 2 linear switching neurons $y=0$}
%\label{fig:t4}
%\end{figure}

\begin{figure}
\begin{center}
\includegraphics[scale=0.3]{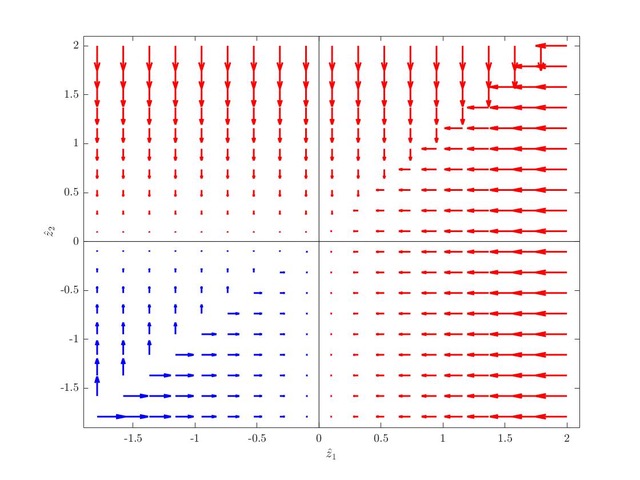}
\end{center}
\caption{Flow field of 2 linear switching neurons $y=1$}
\label{fig:t5}
\end{figure}

\begin{figure}
\begin{center}
\includegraphics[scale=0.3]{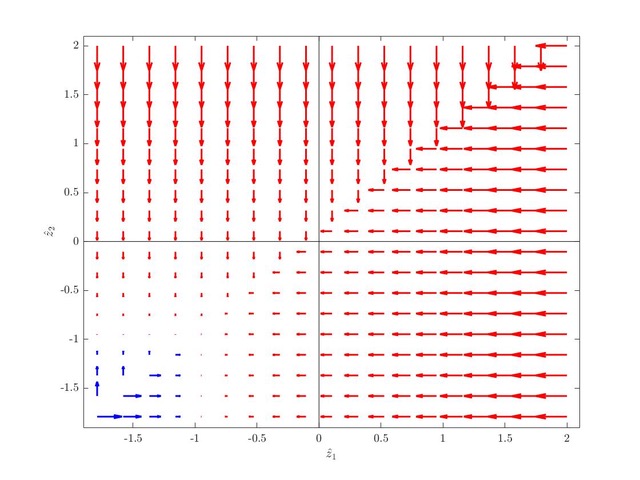}
\end{center}
\caption{Flow field of 2 linear switching neurons $y=-1$}
\label{fig:t6}
\end{figure}

We illustrate this now for the simplified learning law (\ref{eq:dotz}) for $K=1$ for 2 switching linear neurons in figure \ref{fig:t5} and  \ref{fig:t6}. Red coloring means that $\tilde{y}>0$ and blue coloring means that $\tilde{y}<0$. 

All regions are semi-contracting except the $+45$ degree line. Here we define a semi-contracting space, where this semi-diverging half-plane is excluded. All flow lines at the boundary of this semi-contracting space are inflowing which means, that this excluded half-plane will not destabilize the system. More specifically we see that all transitions are semi-contracting for red $\tilde{y}>0$ and semi-diverging for the blue $\tilde{y}<0$. No indifferently frozen regions or separated non-contraction regions exist, i.e. the regions of divergence are a set of zero measure.

Note that the separating half plane between both neurons corresponds to the kink in $\hat{y}_k$ between both linear neurons.

$K>1$ then corresponds to the superposition of several semi-contracting dynamics. This implies that the learning law (\ref{eq:dotb}) is semi-contracting with the switched linear neurons (\ref{eq:lin}).  We summarize this in the following with Theorem \ref{th:theoremF}:
\begin{theorem}
\itshape
The $k$ constant inputs ${\bf x}_k = (1, x_k)^T$ and measured outputs $y_k$ can be approximated with the piece-wise linear function
\begin{equation}
\hat{y}_k =  \hat{z}_{ik} \  \mathrm{: for \  } x_{i-1} \le x \le  x_{i}     \label{eq:linear} \\ 
\end{equation}
with the $i=1, ...,  I$ local neurons  $\hat{z}_{ik} =  \hat{\bf w}_i^T {\bf x}_k$ with $\hat{\bf w}_i = (\hat{w}_{i1}, \hat{w}_{i2})^T$  and ${\bf x}_k = (1, x_k)^T$ between the corners $x_{i-1}, x_{i}$ to the neighboring linear neurons $i-1$ and $i+1$, that solve $\hat{z}_{i} = \hat{z}_{i-1}, \hat{z}_{i} = \hat{z}_{i+1}$ respectively.

The corresponding learning law
\begin{equation}
{\dot{\hat{\bf w}}}_i =  -  \frac{1}{T} \sum_ {k=1}^K \tilde{y}_k \frac{\partial \hat{y}_k}{\partial \hat{z}_{ik}} {\bf x}_k \label{eq:learning1D}
\end{equation}
with $T(t) > 0$ and the measurement error $\tilde{y}_k  = \hat{y}_k - y_k$ is semi-contracting in the invariant set that excludes corners for $\tilde{y}<0$.

In addition contraction is defined in the observability direction of the dyadic block diagonal matrix
\begin{equation}
 - \frac{1}{T} \sum_{k=1}^K  \frac{\partial \hat{y}_k}{\partial \hat{z}_{ik}} \frac{\partial \hat{y}_k}{\partial \hat{z}_{jk}}^T  {\bf x}_k {\bf x}_k^T \le 0 \label{eq:contdir1} 
\end{equation}
 \label{th:learning1D}
 which implies that any trajectory converges to the null-space of $\sum_{k=1}^K  \frac{\partial \hat{y}_k}{\partial \hat{\bf z}_{ik}} \frac{\partial \hat{y}_k}{\partial \hat{\bf z}_{jk}}^T  {\bf x}_k {\bf x}_k^T$. 
  \label{th:1D}
\end{theorem}

Note that separated finited regions do not exist here in contrast to the standard learning of section \ref{standard} since we have $I-1$ separating half-planes in contrast to $I$ separating full-planes at the edges of the approximation function.

Also note the first element in (\ref{eq:contdir1}) defines which neurons are active and the second elements defines if sufficient measurements are given to uniquely define that active linear neuron. This means that the system is contracting if all neurons are active and if sufficient measurements exist to define them locally. Even if the neurons are initialized with full rank in (\ref{eq:contdir1}) it can happen that some neurons become inactive over time, i.e. the neuron range (\ref{eq:linear}) becomes a set of zero measure. This  means the system degrades from contraction to semi-contraction behavior. Such an inactive neuron for all measurements $k$ can be pruned without affecting $\hat{y}_k$ and then contraction behavior is achieved again.

Finally, note that at the same time new neurons can be created around a measurement position ${\bf x}_k = (1, x_k)^T$ with persistent relevant error $\tilde{y}_k  = \hat{y}_k - y_k$, by adding a new corner at $x_k$, i.e. by replacing one linear neuron with two. If we choose $\hat{y}_k$ at this new corner identical to the original function then the cost function (\ref{eq:V}) stays unchanged during neuron creation. This means that the Lyapunov or contraction proof of the cost function (\ref{eq:V}) remains valid during neuron creation since the cost function is unchanged at this time instance. Note that this creation of neurons is necessary to assure that the true weights or function is contained as a special solution in the dynamics. 

\Example{}{Consider the green/blue polygon in figure \ref{fig:ply}, where a measurement $y_k$ is performed at each integer $x_k$. We apply the learning law (\ref{eq:learning1D}) of Theorem  \ref{th:1D} first with a single neuron. This neuron converges to the red line. This red line is then split in 2 lines, after learning (\ref{eq:learning1D}) it converges to the violet line and after a further split to the green line which perfectly matches the blue target function. Here we can see that the neuron splitting or creation is necessary to assure that the true solution is part of the system dynamics.

\begin{figure}
\begin{center}
\includegraphics[scale=0.5]{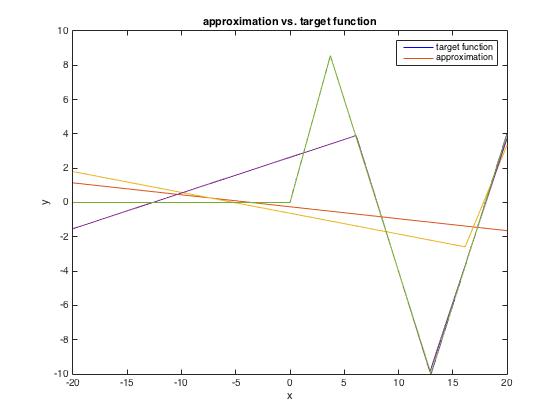}
\end{center}
\caption{Convergence to polygon}
\label{fig:ply}
\end{figure}

}{1dex}

\subsection{Contracting learning of N-dimensional functions} \label{jn}

\begin{figure}
\begin{center}
\includegraphics[scale=0.5]{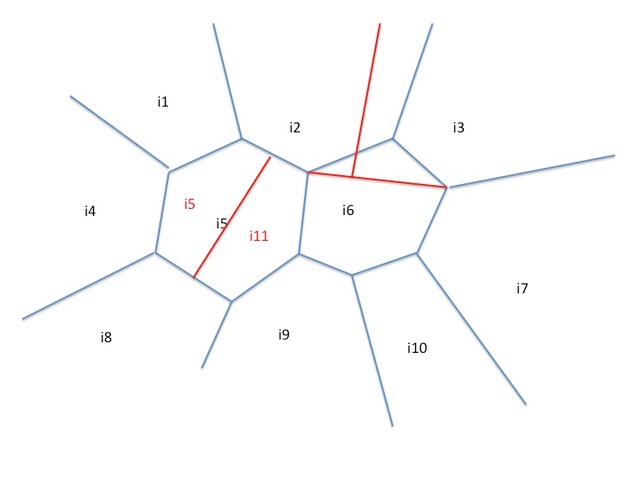}
\end{center}
\caption{Several neurons, where each corner has at most 3 neighboring neurons}
\label{fig:triangle}
\end{figure}

In the following we propose an $N$-dimensional approximation function which is globally semi-contracting except for regions of zero measure. Indeed we locally estimate the piece-wise linear function
\begin{equation}
\hat{y}_k =  \hat{z}_{ik}  \  \mathrm{: within  \ a \ convex \ polyhedron \ } \nonumber
\end{equation}
with the local linear neuron $\hat{z}_{ik} =  \hat{\bf w}_i^T {\bf x}_k$ with $\hat{\bf w}_i = (\hat{w}_{i1}, ..., \hat{w}_{iN+1})^T$ and ${\bf x}_k(t) = (1, x_{1k}^T, .., x_{Nk})^T$. The linear edges of the polyhedron are defined $\forall$ neighbours of the neuron $i$ by
$$
 \hat{z}_{i} = \hat{z}_{neighbor \ of \ i}  
$$
Reformulatilng the above as inequalities allows to check if a point is within or outside a convex polyhedron. As an example figure \ref{fig:triangle} shows in black $i1$ to $i10$ neurons. The neurons, with their neighbours are summarized in table \ref{fig:neig}. 
\begin{table}[htp]
\begin{center}
\begin{tabular}{|c|c|c|c|c|c|c|}
Neuron $i$ & Neighbour of $i$  \\ \hline
i1 & i2 & i5 & i4 & & & \\ \hline
i2 & i3 & i6 & i5 & i1 &  &  \\ \hline
i3 & i7 & i6 & i2 & & & \\ \hline
i4 & i1 & i5 & i8 & & & \\ \hline
i5 & i2 & i6 & i9 & i8 & i4 & i1 \\ \hline
i6 & i3 & i7 & i10 & i9 & i5 & i2 \\ \hline
i7 & i10 & i6 & i3 & & & \\ \hline
i8 & i4 & i5 & i9 & & & \\ \hline
i9 & i8 & i5 & i6 & i10 & & \\ \hline
i10 & i9 & i6 & i7 &  &  &  \\ \hline
\end{tabular}
\end{center}
\caption{Example of neurons with neighbours}
\label{fig:neig}
\end{table}

\begin{figure}
\begin{center}
\includegraphics[scale=0.5]{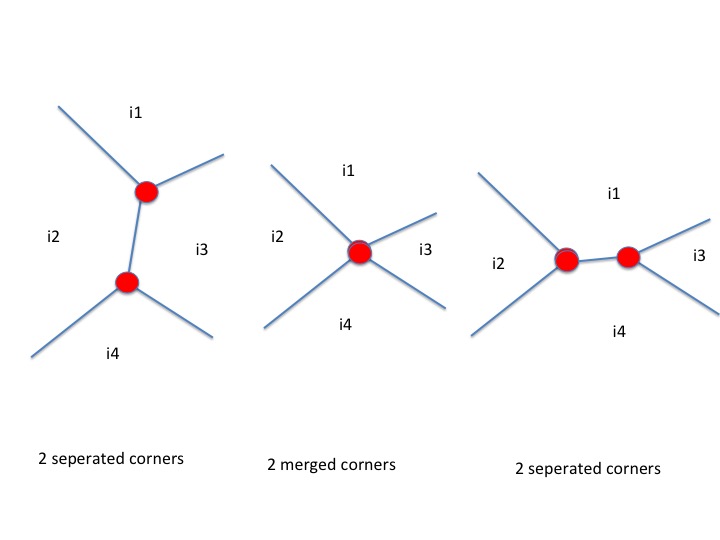}
\end{center}
\caption{Two merged or separated corners}
\label{fig:Seperartion}
\end{figure}

The above leads to the following definition:
\begin{definition}
An $N$-dimensional neuron $i$ is defined by the linear function 
\begin{equation}
\hat{z}_{ik} =  \hat{\bf w}_i^T {\bf x}_k \label{eq:func}
\end{equation}
with $\hat{\bf w}_i = (\hat{w}_{i1}, ..., \hat{w}_{iN+1})^T$ and ${\bf x}_k(t) = (1, x_{1k}^T, .., x_{Nk})^T$ and the neighbouring neurons which define
\begin{itemize}
\item the linear edges of the neuron $i$ to all neighbours 
\begin{equation}
 \hat{z}_{i} = \hat{z}_{neighbor \ of \ i}  \label{eq:edge}
\end{equation}
Reformulating the above as inequality to the interior of the neuron defines the neuron region. 
\item the corner of the neighbouring neurons $i_1, ..., i_{N+1}$, defined by
\begin{equation}
 \hat{z}_{i_1} = .... = \hat{z}_{i_{N+1}} \label{eq:corner}
\end{equation} 
If more than $N+1$ edges cross exactly in one point as in the middle of figure \ref{fig:Seperartion} then a corner can become overdetermined in (\ref{eq:corner}). However, for a limited slope of the neurons this overdetermined corner (\ref{eq:corner}) is solvable in $\hat{y}_k$.

If more than $N+1$ edges happen to cross in several points as in the right or left of figure \ref{fig:Seperartion}, then the neighborhood table may have to be updated.
\end{itemize}
\label{th:neuron}
\end{definition}
We can now compute similarly to the $1$-dimensional case
\begin{equation}
\frac{\partial \hat{y}_k}{\partial \hat{z}_{ik}} = 
\left\{
\begin{array}{l}
1  \mathrm{:for \ active \ neurons \ }  \\
0  \mathrm{:for \ non \ active \ neurons \ }
\end{array} \right. \nonumber
\end{equation}

At the transitions between the linear regions we get the semi-definite Jacobian
\begin{equation}
\frac{\partial^2 \hat{y}_k}{\partial \hat{\bf z}_{ik} \partial \hat{\bf z}_{jk}}  = 
\left\{
\begin{array}{l}
\left(
\begin{array}{ccc}
\delta  &  - \delta  \\
- \delta   &  \delta 
\end{array}
\right) \mathrm{: for \ the  \ corner \ between \ neuron \ } i \ \mathrm{to \ its  \ neighbour } \\
0   \mathrm{: otherwise}  \\
\end{array} \right. \nonumber
\end{equation}
with the Dirac impulse $\delta > 0$. Note that these semi-definite matrices are then weighted with $\tilde{y}_k$.

This implies that the learning law (\ref{eq:dotb}) is semi-contracting.  We summarize this in the following:
\begin{theorem}
\itshape
The $k$ constant inputs ${\bf x}_k(t) = (1, x_{1k}^T, .., x_{Nk})^T$  and measured outputs $y_k$ can be approximated with the piece-wise linear function
\begin{equation}
\hat{y}_k =  \hat{z}_{ik} =  \hat{\bf w}_i^T {\bf x}_k \    \label{eq:lindimJ} \\ 
\end{equation}
with the $i=1, ..., I$ local neurons in the neuron region of Definition \ref{th:neuron}. 

The corresponding learning law

\begin{equation}
{\dot{\hat{\bf w}}}_i =  -  \frac{1}{T} \sum_ {k=1}^K \tilde{y}_k \frac{\partial \hat{y}_k}{\partial \hat{z}_{ik}} {\bf x}_k \label{eq:learningND}
\end{equation}
with $T(t) > 0$ and the measurement error $\tilde{y}_k  = \hat{y}_k - y_k$ is globally semi-contracting except the edges for $\tilde{y}<0$.

In addition, contraction is defined in the observability direction of the dyadic block diagonal matrix
\begin{equation}
 - \frac{1}{T} \sum_{k=1}^K  \frac{\partial \hat{y}_k}{\partial \hat{z}_{ik}} \frac{\partial \hat{y}_k}{\partial \hat{z}_{jk}}^T  {\bf x}_k {\bf x}_k^T \le 0 \label{eq:contdirN} 
\end{equation}
 which implies that any trajectory converges to the null-space of $\sum_{k=1}^n  \frac{\partial \hat{y}_k}{\partial \hat{\bf z}_{k}} \frac{\partial \hat{y}_k}{\partial \hat{\bf z}_{k}}^T  {\bf x}_k {\bf x}_k^T$. 
   \label{th:ND}
\end{theorem}

Note that
\begin{itemize}
\item new neurons can be created around a measurement position ${\bf x}_k$ , with persistent error $\int_0^T | \tilde{y}_k | dt / T > threshold$ with $\tilde{y}_k = \hat{y}_k - y_k $ within a neuron $i$, by splitting an existing neuron. This is shown in figure \ref{fig:triangle} for $i5$ which is split in red $i5$ and $i11$. We choose the neuron function (\ref{eq:lindimJ}) to remain unchanged at splitting such that the Lyapunov or contraction proof of the cost function (\ref{eq:V}) remains valid even during this neuron creation. Both neurons and the neighbours have to update the neighbour-table accordingly.

\item neurons can be pruned if they have the same neuron function (\ref{eq:lindimJ}). This e.g. happens if the edge becomes identical as for the red edge of neuron $i2$ and $i3$ in figure \ref{fig:triangle}. Since both neurons have two identical edges it imples that both could be pruned to one. The red edge between $i2, i3$ and $i6$ is also the transition from a convex to concave region $i6$. The merging of $i2$ and $i3$ avoids that $i6$ can become concave.

In addition neurons can be pruned if the neuron region reduces to a set of zero measure.
\end{itemize}
Following the rules above implies that all neurons stay convex, which simplifies the measurement in/out neuron check. Also the number of neurons is adapted to the needs of the approximated function. 

Note in addition that separated finited regions do not exist here in contrast to the standard learning of section \ref{standard} since we have separating plane parts (local edges) in contrast to seperating full planes (infinite edges) of the approximation function.

\Example{}{As an example figure \ref{fig:corner} shows the convergence of 3 planar neurons to a corner for the learning law (\ref{eq:learningND}) of Theorem  \ref{th:ND}. We can see that the approximation of the 3 neurons to the red corner converged to a vanishing error.

\begin{figure}
\begin{center}
\includegraphics[scale=0.4]{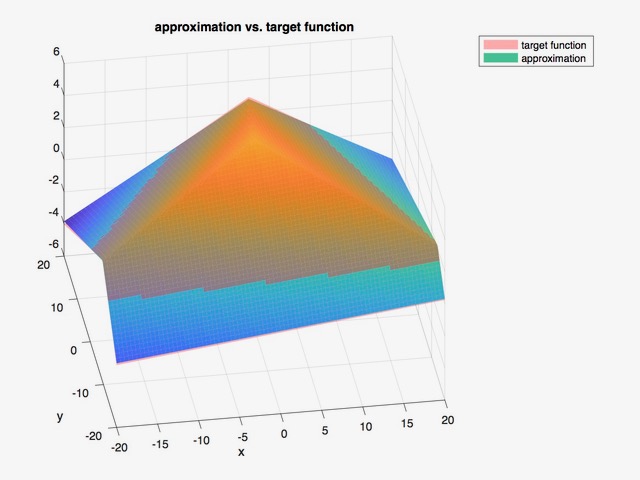}
\end{center}
\caption{Convergence of 3 neurons to a corner}
\label{fig:corner}
\end{figure}

}{Ndex}

\subsection{Summary}

This paper proposes a semi-contracting piece-wise linear neuronal network, where

\begin{itemize}
\item  the neurons in Definition \ref{th:neuron} locally perform 
a semi-contracting linear approximation, 
\item each local approximation occurs in a convex polyhedron, 
\item polyhedron edges are adapted according to the given neighborhood relation
\end{itemize}

\noindent In contrast to standard ReLU-based approaches, learning is
done locally, and explicit stability guarantees are provided.
Creation and pruning of neurons assures that the true approximation is
a special solution of the network. Only a minimal number of neurons
necessary for the approximation is needed.

Extensions to multi-layer networks are the subject of current research.

$ \  $ 

\noindent{\bf Author Contributions}\ \ \ \ WL and JJS developed the theory.
PG performed the simulations.

$  \ $

\end{document}